\documentclass[conference]{IEEEtran} \IEEEoverridecommandlockouts

\usepackage{cite}

\usepackage{graphicx}
\graphicspath{{Images/}}
\DeclareGraphicsExtensions{.pdf,.jpeg,.png,.eps}

\usepackage{amsmath}

\usepackage{algorithmicx}
\usepackage{algpseudocode}
\usepackage{algorithm}
\algnewcommand{\Initialize}[1]{%
  \State \textbf{Initialize:}
  \Statex \hspace*{\algorithmicindent}\parbox[t]{.8\linewidth}{\raggedright #1}
}

\usepackage{array}

\usepackage{stfloats}
\usepackage{url}
\begin{document}
\title{Peer-to-Peer Control of Microgrids}

\author{\IEEEauthorblockN{Hamada Almasalma, Jonas Engels, Geert Deconinck}
\IEEEauthorblockA{ESAT-ELECTA, KU Leuven, Kasteelpark Arenberg 10, 3001 Leuven, Belgium\\
EnergyVille, Thor Park 8300, 3600 Genk, Belgium\\ hamada.almasalma@kuleuven.be, jonas.engels@kuleuven.be}
\IEEEauthorblockA{Proc. IEEE YRS 2016. IEEE Benelux PELS/PES/IAS Young Researchers Symposium. Eindhoven, Netherlands, 12-13 May 2016}
\thanks{This work is partially supported by H2020-LCE-2014-3 Peer to Peer Smart Energy Distribution Networks (P2PSmartTest) project (European Commission).}}


%Telephone: (800) 555--1212, Fax: (888) 555--1212}


\maketitle

\begin{abstract}
%In this paper, the major issues and challenges in control of microgrids are discussed. The motivation to develop microgrids is presented as an effective solution for the control of distribution networks with high level penetration of Distributed Energy Resources (DERs). The paper classifies microgrid control into three levels: primary, secondary and tertiary, %Well, this classification is not a mayor contribution of the paper, we just take this over from literature.where the secondary and tertiary control is of concern of this paper. % do we?
%The paper presents a way to classify the existing secondary and tertiary control techniques from highly centralized into peer-to-peer (P2P) distributed techniques. The paper proposes a P2P voltage control strategy for controlling the microgrid. Different control algorithms suited for the proposed P2P control strategy are discussed and compared.
In this paper, the major challenges and issues in control of microgrids are discussed. The paper classifies possible microgrid control architectures from highly centralized to fully distributed peer-to-peer techniques.
A control paradigm based on coupled microgrids, peer-to-peer communication and autonomous control, is proposed as a way to control the distribution network with a high penetration of distributed energy resources.
The paper suggests epidemic algorithms as an appropriate method for the proposed peer-to-peer control strategy.
\end{abstract}

\vspace{2 mm}
 %there should be space between the Abstract and the keywords?

\section{Introduction} 
Growing concerns about energy sustainability, security of supply and an increasing penetration of renewable energy and other distributed energy resources (DERs), such as storage systems and electrical vehicles, are impacting the operation and the architecture of the electricity system. While currently placing a burden on the distribution grid, it is generally agreed that DERs could also be used for active grid control, thereby contributing to a stable and secure grid. To accomplish this, new control systems have to be designed that are able to fully harness the potential of the installed DERs~\cite{gorman2008enhanced}. This active monitoring and control of the distribution grid is commonly referred to as being essential in the \emph{smart grid}, which is regarded to be key in the future integration of electricity consumers, generators and those that do both (prosumers). %Citation?

Current distribution networks are not designed to accommodate a large amount of DERs. A large penetration of DERs may create problems to maintain the quality of supply to all customers connected to the distribution network. Besides, the intermittent nature of DERs can create issues with the second-by-second balance of demand and supply. However, by coordinating the DERs, these issues could be resolved without the need for additional investments in grid infrastructure~\cite{lopes2007integrating}. 

In the literature~\cite{MG_concept}, the idea of a microgrid is an often mentioned alternative to controlling the whole distribution grid with a large amount of DERs. The main idea is that, when there are many DERs in a wide network, it can be very complex and difficult to control. Thus, a potential way to manage this complexity is by breaking the entire grid down into a smaller microgrids, containing only a limited amount of DERs. This paper elaborates this idea and proposes an operational control paradigm for the future distribution grid, based on the concept of microgrids.

When considering such a microgrid and the coordination of multiple microgrids, different control methods can be found in literature \cite{hatziargyriou2014microgrids}. This paper proposes a taxonomy of these control methods, from fully centralized to completely decentralized control methods. The control architecture that is the most decentralized, the peer-to-peer control architecture, is further elaborated in this paper, as it is a promising way for future control of the distribution grid. Since the DERs are typically highly distributed, operated by many different owners and with different objectives, it is desirable that the control system for a microgrid operates in a highly distributed way as well. Besides, a robust control system is needed that does not depend on a single point of failure, as most of the more centralized control methods do.

%The future smart grid is an electricity network that can intelligently integrate generators, consumers and those that do both (prosumer) in order to efficiently deliver sustainable, economic and secure electricity supplies. % These growing concerns over energy sustainability, security of power supply and the need to increase the shares of renewable energy and Distributed Energy Resources (DERs) are impacting the operation and the architecture of the energy system. ->I put this sentence in front.

%In order to successfully integrate DERs and prosumers, many technical challenges must yet be overcome to ensure that the present levels of reliability are not significantly affected and the potential benefits of DERs and prosumers are fully harnessed. In this sense, the main issues include voltage, frequency and stability issues, which we will elaborate below.

The rest of the paper is organized as follows: section~\ref{Issues} introduces the most prominent issues with regard to the integration of DERs in the electrical grid are, together with a short elaboration of the microgrid concept that would be able to overcome these issues. Section~\ref{Architectures} proposes a categorization for the different architectures of microgrid control. Section~\ref{P2Pparadigm} then proposes a new control paradigm for the distribution grid, based on the microgrid concept and the previously identified peer-to-peer architecture. Finally, the paper is concluded in section~\ref{conclusion}.

%\begin{itemize}

%\item \textbf{Voltage issues:} 
\section{Issues with the Integration of DERs in the Electrical Grid} 
This section summarizes the most prominent issues with regard to the integration of DERs in the current electrical grid. Both voltage issues and frequency or stability issues are discussed. Other type of issues, such as harmonics, security issues and power fluctuations are not discussed here. The microgrid concept is presented as a possible solution that is able to overcome these issues.
\label{Issues}
\subsection{Voltage Issues}
%A large penetration of DERs may create problems to maintain the quality of supply voltage of all customers connected to the same part of the distribution network. Up till now, the distribution grid is designed in such a way that, during moments of minimum load, the voltage received by the last customers on a distribution feeder is just below the maximum allowed. If then e.g. renewable generation is connected at the end of the feeder, the power flow in the circuit will reverse, from the customers to the main grid \cite{Masters2002}. Under such light load conditions, this will lead to the rising of the voltage of the feeder, beyond its allowed limits. The voltage rise in the distribution system limits the hosting capacity of of the DERs. The traditional solution to this problem is to reinforce the local distribution grid by installing more cables. However, this is mostly not cheap, as new infrastructure has to be installed in residential neighborhoods. Another approach is by using the already infrastructure in a more optimal way, by coordination of the local generation, on-load tap changers or other equipment used to control voltage in distribution networks. This would be the purpose of the smart grid The capacity of many traditional distribution networks is limited by the variations in voltage that occur between times of maximum and minimum load and so the circuits are not loaded near to their thermal limits.
A large penetration of DERs may create problems to maintain the voltage quality of all customers connected to the same part of the distribution network. Up till now, voltage quality in the distribution grid is achieved based on the layout of grid infrastructure that is capable of operating within limits even in worst case scenarios, with the assumption of unidirectional power flows. The planning of the infrastructure is quite straightforward: minimum and maximum load conditions are considered and minimum and maximum voltages in the grid are examined. The network is dimensioned in such a way that the minimum voltage is near the lower limit of the allowed voltage range and the maximum voltage is near the upper limit of the allowed voltage range. When connecting significant amounts of distributed generation to the network, the assumption of unidirectional power flows is not always valid any more and the voltage profile of the network can be quite different than in the case without any generation. With maximum load conditions, distributed generation increases the voltage level in the network and, hence, enhance the voltage quality in the grid. However, when the load on the network is at a minimum, the generated power of the distributed generation can reverse the power flows in the grid, what could lead to a rise of the voltage profile beyond its allowed limits. 

The capacity of many traditional distribution networks is limited by the voltage variations that occur between maximum and minimum load conditions. This capacity is dimensioned so that the lines are never loaded beyond their thermal limits. The traditional solution to this problem is to reinforce the local distribution grid by installing more cables. %? hosting capacity or?
However, generally this is quite expensive, as new infrastructure has to be installed in residential neighborhoods. Another approach is by using the already installed infrastructure in a more optimal way, by coordination of the local generation, on-load tap changers or other equipment used to control voltage in distribution networks. This is the purpose of the smart~grid \cite{xu2008voltage}. 

\subsection{Frequency and Stability Issues}
The consumption and generation in the electricity system has to be balanced on a second-by-second basis. Traditionally, this balance is maintained by flexible generation units that are standby and are able to regulate their generated amount of electricity as required. However, the intermittent and unpredictable nature of new renewable energy resources creates issues with this traditional approach, as these new sources are usually not dispatchable. 

Any deviation from the demand and supply balance results in a deviation of the system frequency from its nominal value, while large frequency deviations will affect the system stability. Therefore, the system operators maintain frequency within strict limits by using ancillary services for balancing of the grid, that are able to respond within various time frames \cite{Erinmez1999}. Up till now these balancing services are exclusively organized by the transmission system operator (TSO). However, with the high penetration of DERs, maintaining the supply and demand balance and thus the system frequency within limits becomes more challenging. Spinning reserves or energy storage can address this problem but with a considerable increase in cost. Therefore, power system operators are increasingly seeking new reserves for frequency response from demand flexibility, instead of supply. As most of the demand is connected to the distribution grid, this will mean that the distribution system operators (DSOs) will be involved in this process. 
% \end{itemize}

\subsection{Towards a New Control Paradigm for the Distribution Grid}
One of the key solutions to overcome the above mentioned challenges of integrating DERs in the distribution grid is the design of a smart grid, containing control systems for the coordination of these DERs, thereby ensuring a reliable, secure and economical operation of the distribution network at all times.

Controlling the distribution network to be able to utilize the emerging diversity of DERs at significant levels of penetration, means that the control system has to be able to manage a wide and dynamic set of resources. Controlling such a complex network is not a trivial task. A potential way to manage this complexity is by breaking the entire grid down into a smaller microgrids, containing only a limited amount of DERs. These microgrids should be able to control their local resources as optimal as possible, while being connected to the rest of the grid through a point of common coupling (PCC)~\cite{Lasseter2011}. %Communication with the different DERs needed for controlling them, can than be organised within such a microgrid.

These microgrids can then be coordinated on a higher level into a system of multi-microgrids. The main idea is to design the control of the microgrid in such a way that the microgrid is perceived by the main grid as a single element responding to appropriate control signals. There are different possible architectures of operating such a microgrid, for which this paper proposes a classification ranging from fully centralized methods to fully distributed in section \ref{Architectures}.

A conventional way of controlling such a microgrid, is by a hierarchy of three control levels, each operating on a different timescale and with a different priority: primary, secondary and tertiary~\cite{BidramMicrogrids, Olivares2014}. Primary control is focused on keeping the grid stable in all circumstances, thus it needs the largest priority and should act on the smallest time scale. A robust control method is needed, so that in case e.g. the communication network fails, the primary control is still able to maintain stability of the microgrid. Therefore, as little communication as possible is desired.  %Local control methods seem to be most appropriate to offer primary control services to the grid~\cite{Demirok2011}. -> put this in the section about local control
As primary control is often implemented as some kind of proportional controller, a steady state error remains, that shall have to be eliminated with the use of a secondary control method. Secondary control is activated on a slower time scale, e.g. every 15 minutes, and with a lower priority than primary control. Finally, tertiary control is implemented as the slowest level of control, with as purpose the economically optimal operation of the microgrid. Both secondary and tertiary controls usually require at least some kind of communication, as knowledge about the state of the entire system is needed.

These three control levels can be implemented through various organizational architectures of the microgrid. In the following section, this paper proposes a classification for these different architectures.

 \section{Architectures of Microgrid Control}
 \label{Architectures}
The control of the DERs in a single microgrid can be organised according to many different architectures. They can range from fully centralized control where all decisions are made by a single central controller, to completely decentralized controls where all decisions are made by the local DERs. The required communication architecture changes accordingly. In most practical cases, a hybrid architecture exists, where e.g. primary control is implemented locally, and secondary control centrally.

This paper identifies five different approaches, shown in figure \ref{fig:Architecture}: (a)~centralized control, (b)~hierarchical control, (c)~distributed control, (d)~fully decentralized P2P control and (e)~local control.
\subsection{Centralized Control Architecture}
\begin{figure}[!t]
\centering
\includegraphics[width=\columnwidth]{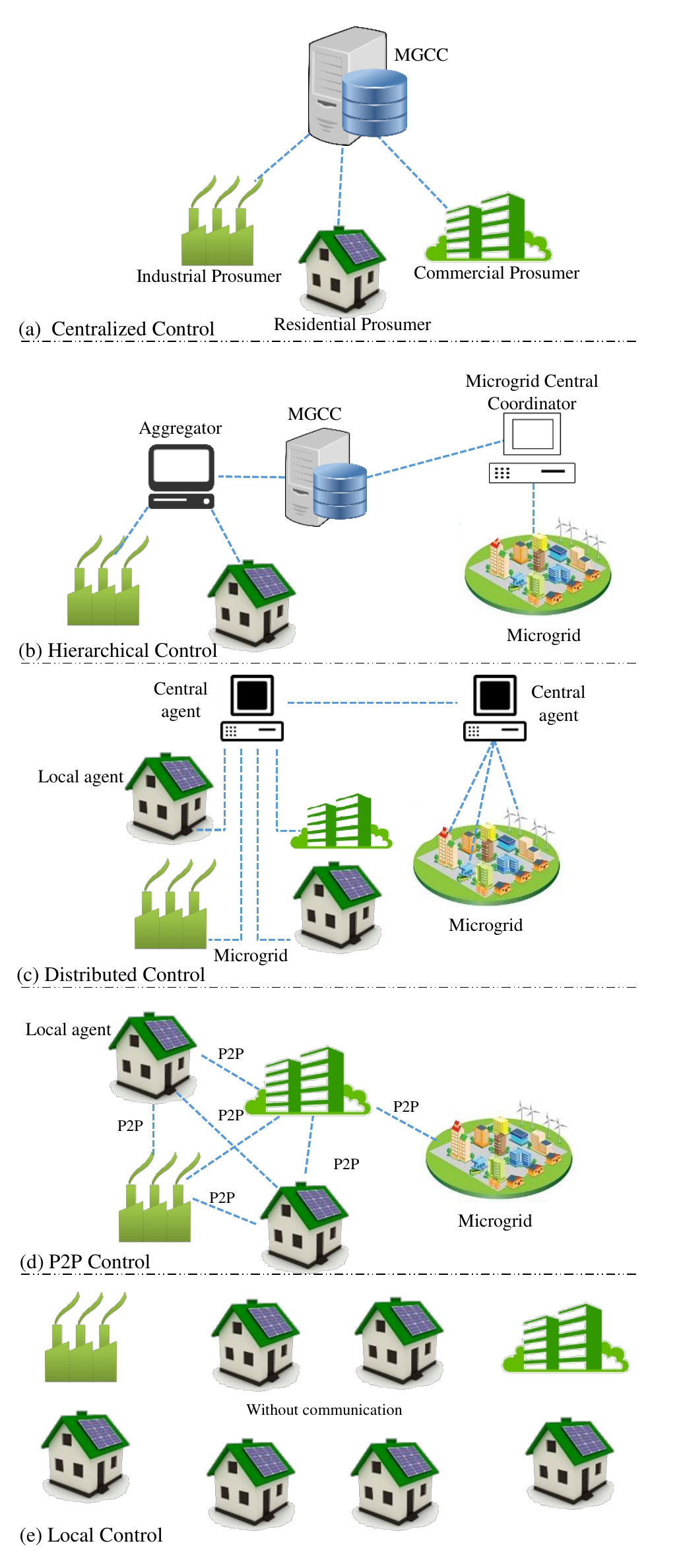}
\caption{Architecture of a microgrids control}
\label{fig:Architecture}
\end{figure}

In an fully centralized design shown in figure~\ref{fig:Architecture}(a), all available measurements %(amplitudes, phase angels and frequency of voltages and currents, active and reactive power of the DERs, ...) 
of the considered microgrid are gathered in a central controller that determines the control actions for all units. Reference~\cite{wamcIP} presents an implementation of a centralized controller based on Wide Area Monitoring and Control system (WAMC) that can be used to implement a centralized secondary and tertiary control. In \cite{wamc_voltage} WAMC has been used to implement a centralized secondary voltage controller. Automatic frequency restoration reserves (aFRR) in the Belgian transmission network, delivered by traditional power plants, are another example of a centrally dispatched  grid control method~\cite{EliaR2}. 

When looking at a single microgrid, the centralized controller is often referred to as a Microgrid Central Controller (MGCC). The advantage of a centralized control system is that the central system receives all necessary data of the microgrid, and based on all available information the multi-objective controller can achieve globally optimal performance. As there is only one controller, this results in a high controllability of the system.
% and high efficiency,
However, this high performance comes at a cost. First of all, the computational burden is heavy, as the optimization is computed based on a large amount of information. Moreover, a centralized controller is a single point of failure and redundancy of the central controller is expensive. The loss of communication with the central controller may cause a shutdown of the overall system. Besides, as all system states and boundary conditions have to be known at the central point, this requires very high quality of communication from all DERs to the central point of control. There is also the concern that the owners of the different DERs are not willing to hand over control of their resources to a third party. Finally, central systems are usually regarded as not being very scalable and system maintenance requires complete shutdown \cite{scc}. To overcome these issues, more distributed control architectures are developed, as described in the sections below.

\subsection{Hierarchical Control Architecture}
%To break the fully centralization topology, the secondary controller can be implemented centrally at the microgrid level as shown in figure~\ref{fig:Architecture}(b). In this case the function of the global controller will be a coordination between the microgrids and the optimization problem. The central controller of each microgrid (MGCC) tries to reach the steady-state values using its local resources and shall ask the global controller for external resources from other microgrids to reach the steady-state values if the internal resources are not sufficient. The use of middle level controllers should increase the reliability of the system. A detailed discussion about the hierarchical centralized topology is presented in [10] and [11].
A first step towards a more decentralized control architecture, is the introduction of a hierarchical system, as shown in figure~\ref{fig:Architecture}(b). Here, there exists some kind of aggregation of the local DERs towards the central controller. Typically, the characteristics of the DERs are represented by a few heuristics or parameters that are combined by an aggregator, who is able to offer these aggregator resources to a central optimizer. The central controller then is able to dispatch the necessary resources, through a hierarchical system of aggregators, that determines which DERs should be used at which moment. Therefore, these methods are also referred to as \emph{Aggregate and Dispatch} methods. As the resources are offered to the central optimizer in an aggregated way, there is considerably less information needed at the central controller, which results in a more scalable system. However, the single point of failure remains, and the points of aggregation might even become new points of failure.

This approach is mostly used for exploiting demand response resources. Examples are \emph{Intelligator}~\cite{Intelligator}, or demand response reserves offered by aggregators to a TSO for e.g. automatic or manual frequency restoration reserves (aFRR or mFRR). The coordination of multiple, centrally controlled microgrids can also be organized in a hierarchical way. In that case, a \emph{microgrid central coordinator} coordinates multiple microgrids, each controlled by a local MGCC. The MGCC of a single microgrid tries to reach an optimal operation point using only its local resources. If the internal resources are not sufficient, the MGCC shall ask the microgrid central coordinator for external resources from other microgrids~\cite{mgcoor}.

%The advantage of these centralized control systems is that the central system receives all the data of system, and then based on the available information the multi-objective controller can achieve global optimal performance and high efficiency. However, this high performance is at high cost. First of all, the computational burden is heavy, as each control signal is computed based on a large amount of information. Moreover, a centralize controller has a single point of failure, a failure in the communication system may cause overall shut down in the system. Another drawbacks of this system are: redundancy of the central controller is expensive, scalability and expansion is a complicated task and system maintenance requires complete shutdown [12]. %-> put this in the centralized control architecture section?

\subsection{Distributed Control Architecture}
In both architectures described above, the DERs are controlled by a third party. However, since the owners of the DERs impose the boundary conditions, one can argue that it might be better to keep the control of the DERs locally. Besides, not all DER owners want to exchange all their information with a third party for privacy reasons. However, to reach a (near) optimal operation of the grid, these DERs should still be coordinated. It is at this point that distributed control architectures come into play. The idea behind distributed control is to divide the centralized problem into a certain number of local controllers or agents. Therefore, each agent does not have a global vision of the problem \cite{distr}, but by means of correct coordination they can reach a globally (near) optimal state.

Coordination is organised by a central agent that is able to communicate global constraints, such as the power limit of a transformer, or exceeding voltage limits. This can be done by the communication of Lagrange multipliers. Examples of algorithms that are suited for this approach are dual decomposition methods or the alternating direction method of multipliers (ADMM)~\cite{DeTurck}. Both are based on the dual ascent method, where price vectors are sent iteratively from the central controller to the DERs. The DERS optimize their consumption towards such a price vector and return demand vectors to the central agent. The central agent then analyses the demand vectors with regard to operational grid constraints, and updates the prices when constraints are being violated. The DERs then optimize again according to this new price vector. This iteration goes on until a steady state solution is found. Figure~\ref{fig:Architecture}(c) represents such a distributed control scheme, consisting of local DERs that optimize and a central agent that controls the global constraints.
%These drawbacks have given motivates to find other approaches that increase the reliability of the control system and reduce its cost [13]. It is at this point where distributed controllers come into play. The idea behind distributed control approaches is to divide the centralized problem to a certain number of local controllers or agents. Therefore, each agent does not have a global vision of the problem [14]. 

Distributed approaches have important advantages that justify their use. As the global optimization problem is divided into several sub-problems, the computational requirements are lower. Besides, the information exchange between local and central agent is limited, which relaxes the requirements of the communication system. This approach results in a very scalable method. As the local DERs perform an optimization by themselves, they do not need to hand over private information to a third party, that controls their resources. However, a central agent still exists, inherently resulting in a single point of failure.

%Their inherent modularity simplifies the system maintenance and the possible expansions of the control system. Moreover, the modularity provides robustness in comparison with a centralized controller. A possible failure does not have to affect the overall system. For this reason, distributed systems have a greater tolerance to failures. Nevertheless, these systems have also some drawbacks that have to be taken into account, being the main one the loss of performance in comparison with a centralized controller. This loss depends on the coordination mechanisms between the agents [15]. The challenge is to design a mechanism that can coordinate the distributed agents in a way that makes the distributed control system to have the same performance and efficiency as the centralized system.

%Figure~\ref{fig:Architecture}(c) presents a distributed control system consists of multi-local agents and a communication agent. The secondary and tertiary control are implemented locally at the agent level. Each agent solves a problem that gives the local references which have to be followed in order to reach the steady state values of voltage and frequency, and the optimized set of the active and the reactive power. Each agent sends values of voltage, frequency, active and reactive power to the communication agent which sends these information to other agents to be able to solve their local problem. This scheme still has a degree of centralization, as losing the communication with the communication agent makes the local agents unable to solve their problems.

\subsection{Peer-to-peer Control Architecture}
\label{sec:p2p}
To eliminate the problems that centralized control inherently possesses, being a single point of failure, the idea of peer-to-peer or autonomous microgrids has been developed. This type of architecture, inspired by P2P computer networking~\cite{P2Pnetworking}, is characterized by the complete absence of a central controller. All local DERs or agents, are equally important and can communicate to other agents~\cite{masdef}, in a peer-to-peer fashion, as shown in figure~\ref{fig:Architecture}(d). The absence of a central controller leads to the term of \emph{autonomous} microgrids. The peer-to-peer communication is used for dissemination of the grid states to all required agents in the microgrid. The grid-supporting agents can then act according to the received information, and in cooperation with each other. In this way they should be able to reach a (near) optimal operation of the considered microgrid. Examples of algorithms that could be used for such P2P communication are gossiping~\cite{KempeGossip} and consensus algorithms~\cite{OlfatiConsensus}. This architecture will be elaborated further in the following sections.

In this architecture, there is a clear absence of a single point of failure. In the case a single agent fails, the other agents can still manage the grid in a stable way. Also when a single communication channel fails, the required information can still reach all necessary participants, via other agents. These properties makes this architecture a robust way of controlling a microgrid. Besides, all information is kept local, eliminating possible privacy concerns. On the other hand, all agents need a considerable amount of local intelligence, as they need to be able to execute the necessary optimizations.

\subsection{Local Control Architecture}
Finally, there also exist control architectures without any form of communication, as shown in figure~\ref{fig:Architecture}(e). This paper classifies them to be \emph{local} control architectures. In this case, optimal operation of the microgrid is rather difficult, as it is impossible to know the complete state of the grid and all operational boundary conditions of the DERs. However, this method is robust against communication failerus, as the absence of any communication will ensure that the grid is also controlled when all communication channels fail. As primary frequency control should be able to operate even when communication fails, this is mostly implemented as a local control method. Thereby it uses droop characteristics that vary active power with variations in locally measured frequency. Another example is local voltage control, implemented by a voltage-reactive power droop~\cite{Demirok2011}.

\section{Proposed Peer-to-peer Based Control Paradigm for the Distribution Grid}
\label{P2Pparadigm}

\begin{figure}[!t]
\centering
\includegraphics[width=\columnwidth]{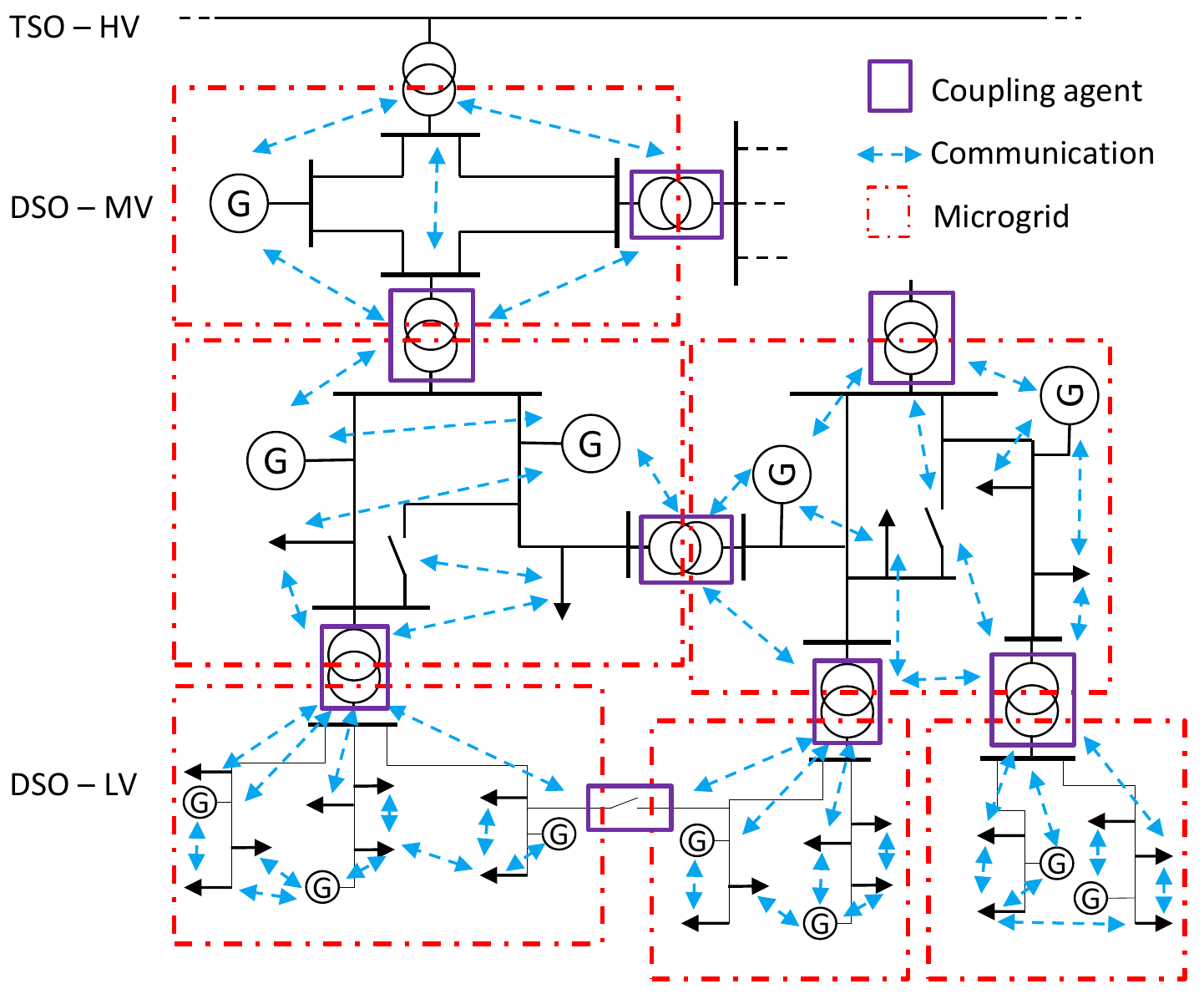}
\caption{P2P based control paradigm}
\label{Controlparadigm}
\end{figure}

As the new DERs are typically highly distributed in the grid, operated by a lot of different owners and with different objectives, it is desirable that the control system for a microgrid is operating in a highly distributed way as well. Plug and play of new resources in this microgrid is a crucial to allow for seamless integration over time. In this perspective, a peer-to-peer control architecture, as introduced in section~\ref{sec:p2p} seems to be a good method for controlling DERs in the distribution grid. It is a robust method and is able to work in a distributed way, without the need for a central controller, having inherent plug and play characteristics. Each agent communicates directly with the agents around it, without having to go through a central device \cite{towards_p2p}. 

Of course, it is impossible to impose this architecture on the whole distribution grid, as it incorporates thousands of DERs that are geographically very dispersed. To deal with this, it is generally agreed that breaking the complete grid down into a smaller microgrids, containing only a limited amount of DERs can be a solution. These microgrids operate according to the presented peer-to-peer control architecture. Points of common coupling are used to connect the different microgrids. 

The proposed scheme is shown in figure~\ref{Controlparadigm}. The distribution network is divided into several microgrids, hierarchically organized on different voltage levels. Each microgrid consists of several autonomous agents. Such an agent could be a renewable generation unit, a group of intelligent controllable loads, a substation, or any other form of DERs. On the connection points of two microgrids there is a coupling agent which serves as gateway of one microgrid to the other microgrid, the point of common coupling. As the microgrids represented in this figure are divided according to their voltage level, a good consideration of such a coupling agent would be the transformer or substation connecting two microgrids. Such a coupling agent represents the characteristics of the whole lower level microgrid (e.g. a low voltage feeder) on the higher level microgrid (e.g. a medium voltage distribution grid). 

The agents are able to communicate with each some neighbouring agents in a peer-to-peer way, creating possibilities to disseminate data about the state of the grid, without the need for one central point of information. Certain algorithms from the field of distributed computing, especially epidemic algorithms, as consensus, gossiping, seem to be particularly appropriate for this goal. They aim at disseminate and aggregate data in distributed networks in a quick and robust way.

\subsection{Epidemic Algorithms Suited for Peer-to-Peer Control}
Epidemic algorithms are used for scalable and efficient data dissemination in distributed P2P networks, without a central controller~\cite{Epidemic}. Examples of such algorithms are consensus and gossiping algorithms, that will be elaborated below:

\subsubsection{Consensus Algorithms}
The main rationale behind consensus algorithms starts from a set of dynamic agents $i$ with an initial state $x_i = z_i$ that communicate this state with some neighbouring agents according to a communication graph. As all agents do this at various times $t$, the states of all agents are spread over the graph in such a way that after convergence, all agents have obtained the same value. With the right properties, this common value will be the average of all initial states of the other agents.

\begin{figure}[!t]
\centering
\includegraphics[scale=0.4]{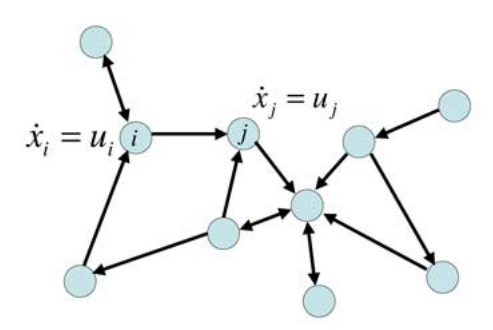}
\caption{Example of a directed graph $G$~\cite{OlfatiConsensus}}
\label{Directedgraph}
\end{figure}

Consider a network of dynamic agents, represented by the directed graph $G=(V,E)$, with the set of nodes $V=\{1,2,…,n\}$ and edges $E\subseteq V \times V$, as depicted in figure~\ref{Directedgraph}. The neighbors of agent $i$ are denoted by $N_i=\{j\in V : (i,j) \in E\}$. The basic consensus algorithm can be described in discrete time as~\cite{OlfatiConsensus,DCConsensus2015}:
\begin{equation}
    x_i[k+1] = x_i[k] + \varepsilon \cdot \sum_{j=1}^n a_{ij} \cdot (x_j[k]-x_i[k]), 
\end{equation}
where $a_{ij}$ indicates the connection status between node $i$ and node $j$, having $a_{ij}=1$ when the nodes not linked and $a_{ij}=0$ when they are not, and $\varepsilon$ the constant edge weight used for tuning of the consensus algorithm. It can be shown that the equilibrium $x^*=(\alpha,\dots,\alpha)$ is globally exponentially stable, with a consensus value $\alpha = 1/n \cdot \sum_i z_i$, being the average of all initial values. Promising research applies these algorithms to DC microgrids in a completely autonomous control way~\cite{DCConsensus2015,DCConsensus2014}. In~\cite{LoadRestoration}, this algorithm is used for global information discovery for automatic load restoration to isolate faults. In \cite{cao2014consensus,Incrementalcost}, a consensus algorithm is used for tertiary control: economic dispatching and equalizing marginal cost over all agents.

\subsubsection{Gossiping Algorithms}
Some epidemic algorithms do not calculate averages of the states of all agents, but rather aggregate values, such as the push-sum gossip-based algorithm\cite{KempeGossip}. Examples of gossip-based control algorithms can be found in~\cite{bolognani2011gossip,debrandere2006voltage}.

\section{Conclusion}
\label{conclusion}
This paper discussed the major issues with regard to the integration of DERs in the electrical grid.
As there are different methods for organizing microgrid control, the paper presented a way to classify 
these methods from highly centralized to completely
distributed P2P techniques. The drawbacks of the centralized control and the advantages of the distributed 
control have been discussed as a motivation to adopt a new control method based on P2P communication and autonomous agents. 
The concept of coupled microgrids is suggested as a way to control the distribution network with penetration of DERs.
Consensus and gossiping algorithms have been proposed as an appropriate techniques for the proposed P2P control paradigm.



%\bibliographystyle{IEEEtran}

\bibliographystyle{IEEEtran}
\bibliography{IEEEabrv,YRSBibliography}

\end{document}